\magnification=\magstep1
\hsize=16truecm
 
\input amstex
\TagsOnRight
\parindent=20pt
\parskip=3pt plus 1.5pt
\define\({\left(}
\define\){\right)}
\define\[{\left[}
\define\]{\right]}
\define\e{\varepsilon}
\define\oo{\omega}

\define\supp {\sup\limits}

\define\summ{\sum\limits}

\define\bigcupp{\bigcup\limits}

\centerline{\bf Sharp estimate on the supremum of a class of} 
\centerline{\bf partial sums of small i.i.d. random variables.}

\smallskip
\centerline{\it P\'eter Major}
\centerline{\it Alfr\'ed R\'enyi Mathematical Institute of the
Hungarian Academy of Science}
\centerline {e-mail address: major.peter$\@$renyi.mta.hu}

\medskip\noindent
{\bf Summary.} {\it We take an $L_1$-dense class of functions 
$\Cal F$ on a measurable space $(X,\Cal X)$ together 
with a sequence of  independent, identically distributed 
$X$-space valued random variables $\xi_1,\dots,\xi_n$
and give a good estimate on the tail distribution of 
$\supp_{f\in\Cal F}\summ_{j=1}^n f(\xi_j)$ if the 
expected values $E|f(\xi_1)|$ are very small for all
$f\in\Cal F$. In a subsequent paper~[2] we shall give 
a sharp bound for the supremum of normalized sums of 
i.i.d. random variables in a more general case. But that 
estimate is a consequence of the results in this work.} 

\beginsection 1. Introduction.

This work is part of a more general investigation 
about the supremum of (normalized) partial sums of 
bounded, independent and identically distributed random 
variables if the class of random variables whose partial 
sums we investigate have some nice properties. It 
turned out that it is useful to investigate first the
case when the expectations of the absolute value of these 
random variables are very small, and this is the subject
of the present paper. In paper~[2] we shall get good 
estimates in the general case when the expectations
of the absolute value of the summands may be relatively 
large with the help of the main result in this paper.

First I recall the notion of $L_1$-dense classes of 
functions which plays an important role in our 
investigation, and then I formulate the main result of 
this paper. After its formulation I make some comments 
that may help in understanding its content and the 
motivation behind this investigation.

\medskip\noindent
{\bf Definition of $L_1$-dense classes of functions.}
{\it Let a measurable space $(X,{\Cal X})$ be given 
together with a class of ${\Cal X}$ measurable, real 
valued functions $\Cal F$ on this space. The class of 
functions ${\Cal F}$ is called an $L_1$-dense class of 
functions with parameter~$D$ and exponent~$L$ if for 
all numbers $0<\varepsilon\le1$ and probability 
measures $\nu$ on the space $(X,{\Cal X})$ there 
exists a finite $\varepsilon$-dense subset
${\Cal F}_{\varepsilon,\nu}=\{f_1,\dots,f_m\}\subset {\Cal F}$
in the space $L_1(X,{\Cal X},\nu)$ with 
$m\le D\varepsilon^{-L}$ elements, i.e. there exists 
such a set ${\Cal F}_{\varepsilon,\nu}\subset {\Cal F}$ 
with $m\le D\varepsilon^{-L}$ elements for which
$\inf\limits_{f_j\in{\Cal F}_{\varepsilon,\nu}}\int |f-f_j|\,d\nu
<\varepsilon$ for all functions $f\in {\Cal F}$. (Here the set
${\Cal F}_{\varepsilon,\nu}$ may depend on the measure $\nu$, 
but its cardinality is bounded by a number depending only on 
$\varepsilon$.)}

\medskip
The main result of this work is the following Theorem~1.

\medskip\noindent
{\bf Theorem 1.} {\it  Let $\Cal F$ be a finite or countable 
$L_1$-dense class of functions with some parameter $D\ge1$ and 
exponent $L\ge1$ on a measurable space $(X,\Cal X)$ such that 
$\sup\limits_{x\in X} |f(x)|\le 1$ for all $f\in\Cal F$. Let 
$\xi_1,\dots,\xi_n$, $n\ge2$, be a sequence of independent and 
identically distributed random variables with values in the 
space $(X,\Cal X)$ with such a distribution $\mu$ for which 
the inequality $\int |f(x)|\mu(\,dx)\le\rho$ holds for all 
$f\in\Cal F$ with a number $0<\rho\le n^{-200}$. Put 
$S_n(f)=S_n(f)(\xi_1,\dots,\xi_n)=\summ_{j=1}^nf(\xi_j)$ for 
all $f\in\Cal F$. The inequality 
$$
P\left(\sup_{f\in \Cal F}|S_n(f)|\ge u\right)\le D\rho^{Cu} \quad
\text{for all } u>41L \tag1.1
$$
holds with some universal constant $1>C>0$. We can choose e.g. 
$C=\frac1{50}$.}

\medskip
I introduce an example that may help in understanding better 
the content of Theorem~1. In particular, it gives some 
hints why a condition of the type $u>CL$ was imposed in 
formula~(1.1). (We applied this condition with $C=41$.)

Let us take a set $X=\{x_1,\dots,x_N\}$ with a large number $N$ 
together with the uniform distribution $\mu$ on it, i.e. let 
$\mu(x_j)=\frac1N$ for all $1\le j\le N$, and define the following
class of function $\Cal F$ on $X$. Fix a positive integer $L$,
and let the class of functions $\Cal F$ consist of the indicator
functions of all subsets of $X$ containing no more than $L$ points.
Let us fix a number~$n$, and choose for all numbers $j=1,\dots,n$ 
a point of the set $X$ choosing each point with the same 
probability $\frac1N$ independently of each other. Let $\xi_j$ 
denote the element of $X$ we chose at the $j$-th time.  In such 
a way we defined a sequence of independent random variables 
$\xi_1,\dots,\xi_n$ on $X$ with distribution $\mu$, and a class 
of functions $\Cal F$ consisting of non-negative functions 
bounded by~1 such that $\int f(x)\mu(\,dx)=\frac LN$ for all 
$f\in\Cal F$. Let us introduce the random sums 
$S_n(f)=\summ_{j=1}^n f(\xi_j)$ for all $f\in\Cal F$. We shall 
estimate first the probability 
$P_n=P\(\supp_{f\in\Cal F}S_n(f)\ge n\)$ and then the probability
$P_{u,n}=P\(\supp_{f\in\Cal F}S_n(f)\ge u\)$ for $u\le n$.

It is not difficult to see that $P_n=1$ if 
$n\le L$, and $P_n\le\binom NL(\frac LN)^n\le C^L\rho^{n-L}$
with $\rho=\frac LN$, where $C$ is a universal constant. The 
number $C$ can be chosen as such a constant for which the 
inequality $p^p\le C^pp!$ holds for all positive integers 
$p$. We can choose for instance $C=4$. In the proof of the 
above estimate we have exploited that $X$ has $\binom NL$ 
subsets containing exactly $L$ points, and the event 
$\supp_{f\in\Cal F}S_n(f)\ge n$ may occur only if there 
is a subset of $X$ with $L$ points such that all $\xi_j$, 
$1\le j\le n$, are contained in this subset. Also the 
estimate $P_{u,n}\le \binom nu P_u\le C^Ln^u\rho^{u-L}$ 
holds, because the event $\supp_{f\in\Cal F}S_n(f)\ge u$ 
can only happen if there are some indices 
$1\le j_1<j_2<\cdots<j_u\le n$ such that all points 
$\xi_{j_s}$, $1\le s\le u$, are contained in a subset 
of $X$ of cardinality $L$. The probability of such an 
event is $P_u$ for all sequences 
$1\le j_1<j_2<\cdots<j_u\le n$, and there are $\binom nu$ 
such sequences.

We show that if $N\ge n^{201}$ and $n\ge41L$, then the
above model satisfies the conditions of Theorem~1, and 
compare the bound we got for $P_{u,n}$ in our previous
calculation with the estimate Theorem~1 supplies in this 
example. To show that the conditions of Theorem~1 hold
in this case we have to prove that the class of functions 
$\Cal F$ consisting of the indicator functions of all 
subsets containing $L$ points of a set $X$ is an 
$L_1$-dense class, and to estimate the probability 
$P_{u,n}$ with the help of Theorem~1 we have to give 
a possible value for the parameter and exponent for this 
$L_1$-dense class. To do this I recall the definition 
of Vapnik--\v{C}ervonenkis classes together with a 
classical result about their properties.

\medskip\noindent
{\bf Definition of Vapnik--\v{C}ervonenkis classes.} 
{\it Let a set $X$ be given, and let us select a class 
${\Cal D}$  of subsets of this set $X$. We call ${\Cal D}$ a 
Vapnik--\v{C}ervonenkis class if there exist two real numbers 
$B$ and $K$ such that for all positive integers $n$ and
subsets $S(n)=\{x_1,\dots,x_n\}\subset X$ of cardinality $n$
of the set $X$ the collection of sets of the form $S(n)\cap D$,
$D\in{\Cal D}$, contains no more than $Bn^K$ subsets of~$S(n)$.
We call $B$ the parameter and $K$ the exponent of this
Vapnik--\v{C}ervonenkis class.}

\medskip
It is not difficult to see that the subsets of a set $X$ 
containing at most $L$ points constitute a 
Vapnik--\v{C}ervonenkis class with exponent $K=L$ and an 
appropriate parameter $B$. (Some calculations show that 
we can choose $B=\frac{1.5}{L!}$.) I would also recall a 
classical result (see e.g.~[3] Chapter 2, 25~Approximation 
Lemma) by which the indicator functions of the sets in a 
Vapnik--\v{C}ervonenkis class constitute an $L_1$-dense 
class of functions. (Actually, the work~[3] uses a slightly 
different terminology, and it presents a more general 
result.) In the book~[3] it is proved that if the parameter 
and exponent of the Vapnik--\v{C}ervonenkis class are $B$ 
and $K$, then the parameter and exponent of the $L_1$-dense 
class consisting of the indicator functions of the sets 
contained in this Vapnik--\v{C}ervonenkis class can be 
chosen as $D=\max(B^2,n_0)$ and $L=2K$ with an appropriate 
constant $n_0=n_0(K)$. But it is not difficult to see by 
slightly modifying the proof that this class of the 
indicator functions can also be considered as an 
$L_1$-dense class of functions with exponent $L=(1+\e)K$ 
and an appropriate parameter $D=D(K,L,\varepsilon)$ for 
arbitrary $\varepsilon>0$.

The above considerations show that the class of functions 
$\Cal F$ considered in the above example is an $L_1$-dense 
class of functions with exponent $2L$ and an appropriate 
parameter $D$. It is even an $L_1$-dense class of functions 
with exponent $(1+\e)L$ and an appropriate parameter $D(\e)$ 
for all $\e>0$. This means in particular that Theorem~1 can 
be applied to estimate the probability $P_{u,n}$ if the
numbers $L$, $N$ and $n$ are appropriately chosen. It is
not difficult to see that both Theorem~1 and our previous
argument provide an estimate of the form 
$P_{u,n}\le\rho^{\alpha u}$ with a universal constant 
$0<\alpha<1$, only the parameter $\alpha$ is different in
these two estimates. (Observe that 
$\rho=\frac LN\ge\int f(x)\mu(\,dx)$ for all $f\in\Cal F$ 
in our example.). To see that we proved such an estimate 
for $P_{u,n}$ which implies the inequality 
$P_{u,n}\le\rho^{\alpha u}$ under the conditions of 
Theorem~1 observe that $\rho^{u-L}\le\rho^{40u/41}$, and 
$n^u\le\rho^{-u/200}$. Moreover, it can be seen 
that if we are not interested in the value of the 
universal parameter $\alpha$, then this estimate is sharp. 
I also remark that in our example we can give a useful 
estimate for $P_{u,n}$ (and not only the trivial bound 
$P_{u,n}\le1$) only in the case $u>L$.

The main content of Theorem~1 is that a similar picture arises
if the supremum of the partial sums defined with the help of an
$L_1$-dense class of functions is considered. Namely, Theorem~1
states that if $\Cal F$ is an $L_1$-dense class of functions 
that satisfies some natural conditions, then there are 
universal constants $0<\alpha<1$, $C_1>1$ and $C_2>0$ such 
that $P\(\supp_{f\in\Cal F}S_n(f)>u\)\le D\rho^{-\alpha u}$ 
if $n\ge C_1L$ and $\rho\le n^{-C_2}$. Here we applied the 
notations of Theorem~1. We also gave an explicit value for 
these universal parameters in Theorem~1, but we did not try 
to find a really good choice. It might be interesting to show 
on the basis of the calculation of the present paper that we can 
choose $C_1=1+\e$ or $\alpha=1-\e$ with arbitrary small $\e>0$ 
if the remaining universal constants are appropriately chosen.
 
As the above considered example shows the estimate of 
Theorem~1 holds only if $u\ge CL$ with a number $C>1$. The
other condition of Theorem~1 by which $\rho\le n^{-C_2}$ 
with a sufficiently large number $C_2>0$ can be weakened. 
Actually this is the topic of paper~[2] which is a 
continuation of the present work. In paper~[2] I shall 
consider such $L_1$-dense classes of functions $\Cal F$ 
for which the parameter $\rho$ considered in Theorem~1 
can be relatively large. On the other hand, in~[2] we 
shall consider only such classes of functions $\Cal F$ 
whose elements have the `normalizing property' 
$\int f(x)\mu(\,dx)=0$ for all $f\in\Cal F$. In the 
present work we did not impose such a normalization 
condition, because in the case $\rho\le n^{-\alpha}$ 
with some $\alpha>1$ the lack of normalization has a 
negligible effect. 

\medskip
Theorem 1 will be proved with the help of Theorem~1A formulated
below. After its formulation I shall explain why Theorem~1A 
can be considered as a very special case of Theorem 1.

\medskip\noindent
{\bf Theorem 1A.} {\it Let $X=\{x_1,\dots,x_N\}$ be a finite 
set of $N$ elements, and let $\Cal X$ be the $\sigma$-algebra 
consisting of all subsets of $X$. Let $\mu$ denote the uniform 
distribution on $X$, i.e. let $\mu(A)=\frac{|A|}N$ for all 
sets $A\subset X$, where $|A|$ denotes the cardinality of a 
set $A$. Let $\Cal F$ be an $L_1$-dense class of functions 
with some parameter $D\ge1$ and exponent $L\ge1$ on the 
measurable space $(X,\Cal X)$ such that 
$0\le f(x)\le 1$ for all $x\in X$ and $f\in\Cal F$, and 
$\int f(x)\mu(\,dx)\le\frac\rho2$ for all $f\in\Cal F$ with 
some $\rho>0$ which satisfies the inequality
$\rho\le\min(\frac1{1000},L^{-20})$. 
Introduce for all numbers $p=1,2,\dots$ the 
$p$-fold direct product $X^p$ of the space $X$ together 
with the $p$-fold  product measure $\mu_p$ of the uniform
 distribution $\mu$ on $X$, i.e. let each sequence 
$x^{(p)}=(x_{s_1},\dots,x_{s_p})$, $x_{s_j}\in X$, 
$1\le j\le p$, have the weight $\mu_p(x^{(p)})=\frac1{N^p}$ 
with respect to the measure $\mu_p$.

For the sake of a simpler argument let us assume that the 
number $N$ has the following special form: $N=2^kN_0$ with 
some integer $k\ge0$, and a number $N_0$ that satisfies the 
inequality $\frac1{16}\rho^{-3/2}<N_0\le\frac18\rho^{-3/2}$. 

Given a function $f\in\Cal F$ and a positive integer~$p$ let 
us define the set $B_p(f)\subset X^p$ for all $p\ge2$ by 
the formula 
$$
B_p(f)=\{x^{(p)}=(x_{s_1},\dots,x_{s_p})\colon\; x^{(p)}\in X^p,
\quad f(x_{s_j})=1 \quad\text{for all }1\le j\le p\}, \tag1.2
$$
and put
$$
B_p=B_p(\Cal F)=\bigcup_{f\in\Cal F} B_p(f). \tag1.3
$$
If $p\ge 2L$ and $p\le\rho^{-1/100}$, then there exist some 
universal constants $C_1>0$ and $1>C_2>0$ such that 
$$
\mu_p(B_p)=\mu_p(B_p(\Cal F))\le C_1D\rho^{C_2p}. \tag1.4
$$
We can choose for instance $C_1=2$ and $C_2=\frac14$.}

\medskip
In Theorem~1A we considered a very special case of the 
problem discussed in Theorem 1. We took a space of the 
form $X=\{x_1,\dots,x_N\}$ with the uniform distribution 
$\mu$ on it, and considered an $L_1$-dense class of functions 
with some special properties. If we apply it with the choice 
$p=n$, then the event $B_p(\Cal F)$ defined in~(1.3) agrees 
with the event $\supp_{f\in\Cal F}S_n(f)\ge n$, and formula~(1.4) 
implies the estimate~(1.1) with the special choice $u=n$ for 
the system $X$, $\Cal F$, $\mu$ considered in Theorem~1A.

Theorem 1A can be proved by means an appropriate induction,
where we can exploit the $L_1$-dense property of the class
of functions $\Cal F$. This will be done in Section~2. In
Section~3 we prove Theorem~1 with the help of Theorem~1A
and a good approximation.

\beginsection 2. The proof of Theorem 1A.

Theorem~1A will be proved by means of induction with respect
to the parameter~$k$ (appearing in the definition of the 
size~$N$ of the set~$X$). The first result of this section, 
Lemma~2.1, formulates a result similar to Theorem~1A in the 
special case when the set $X$, where the functions $f$ are 
defined contains relatively few points. We need it to start 
our induction procedure. 

\medskip\noindent
{\bf Lemma 2.1.} {\it Let us fix a number $\rho$, $0<\rho<1$,
and a set $X=\{x_1,\dots,x_{N_0}\}$, with 
$N_0\le\frac18\rho^{-3/2}$ points together with a class of
functions $\Cal F$ defined on $X$ which satisfies the 
following weakened version of the $L_1$-dense property with 
parameter~$D\ge1$ and exponent~$L\ge1$. For all $0\le u\le1$ 
there is a set of functions $f_1,\dots,f_s$ from the class of 
functions $\Cal F$ with $s\le Du^{-L}$ elements in such a way
that $\inf\limits_{1\le j\le s}\int |f-f_j|\,d\mu\le u$, where
$\mu$ denotes the uniform distribution on $X$. Let us also
assume that $\int f(x)\,d\mu(x)\le\rho$ and $f(x)\ge0$ for
all $f\in F$ and $x\in X$. Let us consider an integer $p\ge 2L$, 
the set $B_p=B_p(\Cal F)\subset X^p$ introduced in formula~(1.3) 
together with the uniform measure $\mu_p$ on the $p$-fold
product $X^p$ of the space~$X$. The inequality
$$
\mu_p(B_p)\le D\rho^{p/4} \tag2.1
$$
holds.}

\medskip\noindent
{\it Proof of Lemma~2.1.}
Let us choose such a set of functions $f_1,\dots,f_s$, 
$f_j\in\Cal F$ for all $1\le j\le s$, with cardinality 
$s\le D\cdot(2N_0)^L$, which has the property that for all 
$f\in\Cal F$ there is a function $f_j$, $1\le j\le s$, for which 
the inequality $\int |f(x)-f_j(x)|\mu(\,dx)\le \frac1{2N_0}$
holds. If $\int |f(x)-f_j(x)|\mu(\,dx)\le \frac1{2N_0}$, then 
$|f(x)-f_j(x)|\le \frac12$ for all $x\in X$. This follows
from the inequality $\frac1{N_0}|f(x)-f_j(x)|\le
\int |f(x)-f_j(x)|\mu(\,dx)\le \frac1{2N_0}$ for all $x\in X$. 
As a consequence,
$\{x\colon\; f(x)=1\}\subset \{x\colon\; f_j(x)\ge\frac12\}$
for such a pair of functions $f$ and $f_j$, and 
$$
B_p=B_p(\Cal F)=\bigcup_{f\in\Cal F} B_p(f)
\subset \bigcup_{j=1}^s\left\{(x_{t_1},\dots,x_{t_p})
\colon f_j(x_{t_k})\ge\frac12
\quad \text{for all } 1\le k\le p\right\}.
$$
Besides, we have for each $j$, $1\le j\le s$,
$$
\mu_p\left\{(x_{t_1},\dots,x_{t_p})\colon f_j(x_{t_k})\ge\frac12
\quad\text{for all }1\le k\le p\right\}
=\(\mu\left\{x_t\colon f_j(x_t)\ge\frac12\right\}\)^p\le (2\rho)^p.
$$
Hence the relations $p\ge2L$ and $N_0\le\frac18\rho^{-3/2}$ imply
that 
$$
\mu_p(B_p)\le s (2\rho)^p\le D(2N_0)^{p/2}(2\rho)^p\le D\rho^{p/4}.
$$
Lemma 2.1 is proved.

\medskip
In our inductive proof we also need a result presented in 
Lemma~2.2.  It is a version of the following heuristic 
statement. Let us consider the supremum of the integrals 
$\int f(x)\mu(\,dx)$ for all functions $f\in\Cal F$ of an 
$L_1$-dense class $\Cal F$ of non-negative functions 
bounded by~1 on a finite set $X$ with respect to the 
uniform distribution $\mu$ on $X$. Let the cardinality
of the set~$X$ be $2N$, where the number $N$ is of the 
form $N=A2^k$ with some positive integers $A$ and $k$, 
and let the above supremum of integrals be bounded by a 
number $\rho_{k+1}$. Then there is a number $\rho_k$ 
slightly larger than $\rho_{k+1}$ with the following 
property. For most subsets $Y\subset X$ with cardinality 
$N$ the supremum of the integrals of the restrictions of 
the functions $f\in\Cal F$ to the set $Y$ with respect to 
the uniform distribution on $Y$ can be bounded by $\rho_k$. 

\medskip\noindent
{\bf Lemma~2.2.} {\it Let us define two sequences of 
numbers 
$$
N_k=2^kN_0,\quad \text{and} \quad 
\rho_k=\rho\prod_{j=0}^{k-1}\(1+\frac3{N_j^{1/8}}\)^{-1},
\qquad k=1,2,\dots,\quad \rho_0=\rho, \tag2.2
$$ 
with the help of some starting numbers $N_0$ and $\rho$ which 
satisfy the relations $\rho\le\min(\frac1{1000},L^{-20})$ 
and $\frac1{16}\rho^{-3/2}<N_0\le\frac18\rho^{-3/2}$. 
Let us fix an integer $k\ge0$, and consider a set 
$X=\{x_1,\dots,x_{2N_k}\}$ with $N_{k+1}=2N_k=N_02^{k+1}$
elements together with an $L_1$-dense class of functions 
$\Cal F$ on $X$ with parameter $D\ge1$ and exponent 
$L\ge1$ such that $0\le f(x)\le 1$ for all points $x\in X$ and 
functions $f\in\Cal F$. Put 
$R_{k+1}(f)=\summ_{j=1}^{N_{k+1}}f(x_j)$, 
and assume that the class of functions $\Cal F$ also 
satisfies the condition $R_{k+1}(f)\le N_{k+1}\rho_{k+1}$ 
for all $f\in\Cal F$. Let us define the quantity 
$R_Y(f)=\summ_{x_j\in Y} f(x_j)$ for all functions 
$f\in\Cal F$ and sets $Y\subset X$. The following 
Statement~(a) holds.

\medskip
\item{(a)} The number of sets $Y\subset X$ such that $|Y|=N_k$, 
and $\supp_{f\in\Cal F} R_Y(f)\ge N_k\rho_k$ is less than 
$\binom{2N_k}{N_k}D\exp\left\{-\frac1{100}2^{k/20}\rho^{-1/20}\right\}$.

}\medskip\noindent
{\it Proof of lemma~2.2.}\/ Let us fix a partition of 
$X=\{x_1,\dots,x_{2N_k}\}$ to two point subsets 
$\{x_{j_1},x_{j_2}\}$,\dots, $\{x_{j_{2N_k-1}},x_{j_{2N_k}}\}$
together with a sequence of iid. random variables 
$\e_1,\dots,\e_{N_k}$ with distribution 
$P(\e_l=1)=P(\e_l=-1)=\frac12$ for all $1\le l\le N_k$. Let 
us define with their help the `randomized sum' 
$$
U_k(f)=\sum_{l=1}^{N_k}\e_l\(f(x_{j_{2l-1}})-f(x_{j_{2l}})\) \tag2.3
$$
for all $f\in\Cal F$. 

Let us observe that for all $f\in\Cal F$ the inequality
$$
P(U_k(f)>2z)\le \exp\left\{-\frac{2z^2}
{\summ_{l=1}^{N_k}(f(x_{j_{2l-1}})-f(x_{j_{2l}})^2}\right\}
\le e^{-z^2/2N_k\rho_{k+1}} 
\quad\text{for all }z>0 \tag2.4
$$
holds by the Hoeffding inequality (see e.g.~[3] Appendix~B)
and the inequality
$$
\sum_{l=1}^{N_k}(f(x_{j_{2l-1}})-f(x_{j_{2l}}))^2
\le2\sum_{j=1}^{2N_k}f(x_j)^2\le 2R_{k+1}(f)\le 4N_k\rho_{k+1}. \tag2.5
$$
(In formula~(2.5) we exploit the condition $0\le f(x)\le1$ which implies
that $f(x_j)^2\le f(x_j)$.)

Define the (random) set 
$V_k=V_k(\e_1,\dots,\e_{N_k})=\bigcupp_{l\colon\;\e_l=1}\{x_{j_{2l-1}}\}
\cup\bigcupp_{l\colon\;\e_l=-1}\{x_{j_{2l}}\}$. With such a 
notation we can write 
$$
\align
\left\{\oo\colon\!\!\!\!\!\!
\summ_{s\in V_k(\e_1(\oo),\dots,\e_{N_k}(\oo))}
\!\!\!\!\!\!\!\!\!\!\! f(x_s) >N_k\rho_{k+1}+z\right\}&\subset
\left\{\oo\colon\!\!\!\!\!\!
\summ_{s\in V_k(\e_1(\oo),\dots,\e_{N_k}(\oo))}
\!\!\!\! \!\!\!\!\!\!\! f(x_s)>\frac{R_{k+1}(f)}2+z\right\} \\ 
&=\{\oo\colon\; U_k(f)(\oo)>2z\}.
\endalign
$$ 
Hence 
$$
P\(\left\{\oo\colon\;\summ_{s\in V_k(\e_1(\oo),\dots,\e_{N_k}(\oo))}
f(x_s)>N_k\rho_{k+1}+z\right\}\)
\le e^{-z^2/2N_k\rho_{k+1}} \quad\text{for all }z>0 \tag2.6
$$
by relation (2.4).

I claim that relation (2.6) implies the following Statement~(b).

\medskip
\item{(b)} For all $f\in\Cal F$ and $z>0$ the number of sets 
$V\subset X$ such that $|V|=N_k$, and 
$\summ_{x\in V} f(x)\ge N_k\rho_{k+1}+z$ is less 
than or equal to $e^{-z^2/2N_k\rho_{k+1}}\binom{2N_k}{N_k}$.

\medskip
Indeed, it follows from relation~(2.6) that for a fixed partition
of the set~$X$ to two point subsets the number of those 
subsets $V\subset X$ which contain exactly one point 
from each element of this partition, (and as a consequence
contain exactly $N_k$ points), and 
$\summ_{s\in V}f(x_s)>N_k\rho_{k+1}+z$ is less than or 
equal to $2^{N_k}e^{-z^2/2N_k\rho_{k+1}}$. We get an
upper bound for the quantity considered in statement~(b) by 
summing up the number of sets $V$ with these properties for 
all partitions of $X$ to two point subsets, and taking into account 
how many times we counted each set $V$ in this procedure. The 
number of the partitions of $X$ to two point subsets equals 
$(2N_k-1)(2N_k-3)\cdots3\cdot1=\frac{(2N_k)!}{2^{N_k}N_k!}$, 
and each partition provides at most 
$2^{N_k}e^{-z^2/2N\rho_{k+1}}$ sets $V$ with the desired 
properties. All sets $V$  were counted $N_k!$-times in this 
calculation. (A set $V$, $|V|=N_k$, was counted in the above 
calculation as many times as the number of those partitions 
of $X$ to two point subsets which have the property that all 
of their elements contain a fixed element of $V$.) These 
considerations imply Statement~(b).

Given a number $0\le u<1$ there exist $s\le Du^{-L}$ functions 
$f_1,\dots,f_s$ in $\Cal F$ with the property that for all 
$f\in\Cal F$ and sets $Y\subset X$ one of the functions $f_j$, 
$1\le j\le s$, satisfies the inequality $\summ_{x\in Y}
|f_j(x)-f(x)|\le\summ_{x\in X}|f_j(x)-f(x)|\le uN_{k+1}$.
We get this relation by applying the $L_1$-density property 
of the class $\Cal F$ (with parameter~$D$ and exponent~$L$) 
with the uniform distribution $\mu$ on $X$. 
This has the consequence that if 
$\summ_{x\in Y}f(x)\ge N_k\rho_{k+1}+z+2uN_k$ for some 
$Y\subset X$ and $f\in \Cal F$, then there exists some 
index $1\le j\le s$ such that
$\summ_{x\in Y}f_j(x)\ge N_k\rho_{k+1}+z$ with the same set 
$Y\subset X$. Hence Statement~(b) implies that the number of 
sets $Y$ such that $|Y|=N_k$ and 
$\summ_{x\in Y}f(x)\ge N_k\rho_{k+1}+z+2uN_k$ with some 
$f\in\Cal F$ is less than or equal to
$s\cdot e^{-z^2/2N_k\rho_{k+1}}\binom{2N_k}{N_k}
=Du^{-L} e^{-z^2/2N_k\rho_{k+1}}\binom{2N_k}{N_k}$.

Put $z=N_k\rho_{k+1}\cdot N_k^{-1/8}$ and $u=\frac z{N_k}$. 
With such a choice we get that the number of sets 
$Y\subset X$ such that $|Y|=N_k$ and 
$\supp_{f\in\Cal F} R_Y(f)\ge 
N_k\rho_{k+1}(1+3N_k^{-1/8})=N_k\rho_k$ 
is less than 
$$
D\(\frac{N_k^{1/8}}{\rho_{k+1}}\)^L
e^{-N_k^{3/4}\rho_{k+1}/2}\binom{2N_k}{N_k}=
\binom{2N_k}{N_k}D\(\frac{2^{k/8}N_0^{1/8}}{\rho_{k+1}}\)^L
e^{-2^{3k/4}N_0^{3/4}\rho_{k+1}/2}. \tag2.7
$$
It follows from the definition of $\rho_k$ that 
$\frac12\rho\le\rho_{k+1}\le\rho$, and we also have 
$L\le \rho^{-1/20}$ because of the condition imposed on the 
number $\rho$. These relations together with the condition 
$\frac1{16}\rho^{-3/2}<N_0\le\frac18\rho^{-3/2}$ of Lemma~2.2
enable us to bound the expression in~(2.7) from above by
$$
\binom{2N_k}{N_k}D
\left(C_12^{k/8}\rho^{-19/16}\right)^{\rho^{-1/20}}
e^{-C_22^{3k/4}\rho^{-1/8}}
\le \binom{2N_k}{N_k}D\exp\left\{-C_32^{k/20}\rho^{-1/20}\right\}
$$
with appropriate constants $C_1$, $C_2$ and $C_3$. One can choose
e.g. $C_3=\frac1{100}$, and this implies Statement~(a). (In the
estimate of the last step we exploited that for a small number 
$\rho>0$ and all positive integers~$k$ the term 
$e^{-C_22^{3k/4}\rho^{-1/8}}$ is much smaller than the
reciprocal of 
$\left(C_12^{k/8}\rho^{-19/16}\right)^{\rho^{-1/20}}$ which is
of order 
$\exp\left\{-\text{const.}\,\rho^{-1/20}(k+\log\frac1\rho)\right\}$.) 
Lemma~2.2 is proved.

\medskip\noindent
{\it Remark.} It may be worth remarking that the most important 
part of Lemma~2.2, relation~(2.4) or its consequence~(2.6) can 
be considered as a weakened version of Lemma~3 in~[1], and even 
its proof is based on the ideas worked out in~[1]. In 
formula~(2.4) a random sum denoted by $U_k(f)$ was estimated by 
means of the Hoeffding inequality. To get this estimate we had 
to bound the variance of the random variable $U_k(f)$, and this 
was done in formula~(2.5). In Lemma~3 of~[1] a similar random 
sum was investigated, but in that case a good asymptotic formula 
and not only an upper bound was proved for the tail distribution 
of the random sum. In the proof of that result a sharp version 
of the central limit theorem was applied instead of the 
Hoeffding inequality, and we needed a good asymptotic formula 
and not only a good upper bound for the variance of the random 
sum we investigated. The proof of the good asymptotic formula 
for this variance was the hardest part in the proof of Lemma~3 
of~[1].

\medskip\noindent
{\it Proof of Theorem~1A.}\/ Let us fix some numbers $N_0$, 
$\rho$ and $L$ which satisfy the conditions of Lemma~2.2. 
Take an integer $k\ge0$, define the numbers $N_k$ and $\rho_k$ 
by formula~(2.2), consider a space  $X=\{x_1,\dots,x_{N_k}\}$ 
with $N_k$ elements, and an $L_1$-dense class of functions 
$\Cal F$ on it with parameter $D\ge1$ and exponent $L\ge1$ 
such that $0\le f(x)\le1$ for all $x\in X$ and $f\in\Cal F$, 
and $\int f(x)\mu(\,dx)\le \rho_k$ for all $f\in\Cal F$ with 
the uniform distribution~$\mu$ on~$X$. Fix an integer $p$ 
such that $p\ge2L$, $p\le\rho^{-1/100}$, and let us also 
consider the sets $B_p(f)$, $f\in\Cal F$, and 
$B_p=B_p(\Cal F)$ introduced in formulas~(1.2) and~(1.3). 
They consist of sequences 
$x^{(p)}=(x_{s_1},\dots,x_{s_p})\in X^p$ with some nice 
properties. Let $V(p,\rho,N_0,k)=V_{D,L}(p,\rho,N_0,k)$ 
denote the supremum of the cardinality of the sets 
$B_p(\Cal F)$ if the supremum is taken for all possible 
sets $X$ and class of functions $\Cal F$ with the above 
properties (with parameters $N_k$ and $\rho_k$).

I claim that
$$
V(p,\rho,N_0,k)\le C_k N_k^p D\rho^{p/4} \quad \text{for all }
k=0,1,2,\dots \tag2.8
$$
with 
$$
C_k=\prod_{j=0}^k (1+2^{-j}\rho). \tag2.9
$$

Relation (2.8) will be proved by means of induction with respect to
$k$. Its validity for $k=0$ follows from Lemma~2.1. Let us assume
that it holds for some $k$, take a set $X$ with cardinality 
$N_{k+1}=2N_k$ together with a class of functions $\Cal F$ which
satisfies the above conditions with the parameters $D$, $L$, $p$,
$\rho_{k+1}$ and $N_{k+1}$, and let us give a good bound on the 
cardinality of the set $B_p(\Cal F)$ defined in~(1.2) and~(1.3) 
in this case. To calculate the number of sequences 
$x^{(p)}=(x_{s_1},\dots,x_{s_p})\in X^p$ which belong to the set 
$B_p(\Cal F)$ let us take all sets $Y\subset X$ with cardinality 
$|Y|=N_k$, let us bound the number of those sequences 
$x^{(p)}\in B_p(\Cal F)$ for which also the property 
$x^{(p)}\in Y^p$ holds, and let us sum up these numbers for 
all sets $Y\subset X$ such that $|Y|=N_k$. Then take into account 
how many times we counted a sequence $x^{(p)}$ in this summation. 
I claim that we get the following estimate in such a way:
$$
|B_p(\Cal F)|\le N_k^p \frac{\binom{2N_k}{N_k}}{\binom {2N_k-p}{N_k-p}}
\left(C_k D\rho^{p/4}
+D\exp\left\{-\frac1{100}2^{k/20}\rho^{-1/20}\right\}\right) \tag2.10
$$
with the coefficient~$C_k$ defined in~(2.9).

To prove relation~(2.10) let us first observe that if $\Cal F$ 
is an $L_1$-dense class of functions on the set $X$ with 
parameter~$D$ and exponent~$L$, and we restrict the domain 
where the functions of $\Cal F$ are defined to a smaller set 
$Y\subset X$ then the class of functions we obtain in such a 
way remains $L_1$-dense with the same parameter~$D$ and 
exponent~$L$. Hence if we fix a set $Y$ with cardinality 
$|Y|=N_k$ for which the property 
$\supp_{f\in\Cal F}R_Y(f)\le N_k\rho_k$ holds (with the 
quantity $R_Y(f)$  introduced in the formulation of 
Lemma~2.2), then the number of those sequences $x^{(p)}$ 
for which $x^{(p)}\in B_p(\Cal F)\cap Y^p$ can be bounded 
by our induction hypothesis by $C_k N_k^p D\rho^{p/4}$. 
We shall bound the number of the sequences 
$x^{(p)}\in\Cal B_p(\Cal F)\cap Y^p$ for the remaining 
sets $Y$ with cardinality $|Y|=N_k$ by the trivial upper 
bound $N_k^p$, but the number of such sets $Y$ is less than
$\binom{2N_k}{N_k}D\exp\left\{-\frac1{100}2^{k/20}\rho^{-1/20}\right\}$
by Lemma~2.2. This yields the upper bound
$C_k N_k^p D\rho^{p/4}\binom{2N_k}{N_k}+
N_k^p\binom{2N_k}{N_k}D\exp\left\{-\frac1{100}2^{k/20}\rho^{-1/20}\right\}$
for the sum we get by summing up the number of sequences 
$x^{(p)}\in Y^p\cap B_p(\Cal F)$ for all subsets with 
$|Y|=N_k$ elements. To prove~(2.10) we still have to take into 
account how many times we counted the sequences 
$x^{(p)}\in B_p(\Cal F)$ in this summation. If all coordinates 
of a sequence $x^{(p)}\in B_p(\Cal F)$ are different, then we 
counted it $\binom{2N_k-p}{N_k-p}$-times, because to find a 
set $Y$, $|Y|=N_k$, containing the elements of this sequence 
$x^{(p)}$ we have to extend these points with $N_k-p$ new 
points from the remaining $2N_k-p$ points of~$X$. If some 
coordinates of a sequence $x^{(p)}$ may agree, then we 
might have counted this sequence with greater multiplicity. 
The above considerations imply~(2.10).

To prove relation~(2.8) with the help of~(2.10) let us observe that
under the conditions of Theorem~1A (In particular, we have 
$\frac1{N_0}\le 16\rho^{3/2}$, 
$p^2\le\rho^{-1/50}\le\frac1{16}\rho^{-1/6}$, $2N_k-p\ge N_k=2^kN_0$ for
all $k=0,1,2,\dots$, and $\rho>0$ is sufficiently small.)
$$
\align
N_k^p \frac{\binom{2N_k}{N_k}}{\binom {2N_k-p}{N_k-p}}
&=N_k^p \frac{\binom{2N_k}{N_k}}{\binom {2N_k-p}{N_k}}
=N_k^p\frac{2N_k(2N_k-1)\cdots (2N_k-p+1)}{N_k(N_k-1)\cdots (N_k-p+1)}\\
&=N_{k+1}^p\left(1+\frac1{2(N_k-1)}\right)\left(1+\frac2{2(N_k-2)}\right)
\cdots\left(1+\frac{p-1}{2(N_k-p+1)}\right)\\
&\le N^p_{k+1}\exp\left\{\frac{p^2}{2^{k+1}N_0}\right\}\le
N^p_{k+1}e^{2^{-(k+1)}\rho^{4/3}}
\le N^p_{k+1}\(1+\frac132^{-(k+1)}\rho\),
\endalign
$$
and
$$
\align
\exp\left\{-\frac1{100}2^{k/20}\rho^{-1/20}\right\}
&=\rho^{p/4} \exp\left\{-\frac1{100}2^{k/20}\rho^{-1/20}
+\frac p4\log\frac1\rho\right\}\\
&\le C_k\rho^{p/4}\cdot\frac13 2^{-(k+1)}\rho
\endalign
$$
with the coefficient~$C_k$ defined in~(2.9). These estimates 
together with (2.10) imply (2.8) for parameter~$k+1$.

It is not difficult to prove Theorem~1A with the help of
relation~(2.8). To do this let us observe that 
$\rho_k\ge\frac\rho2$ and $C_k\le2$ for all $k=0,1,2,\dots$.
Hence taking  a class of functions $\Cal F$ on a set
$X$ with cardinality $N_k$ with some $k\ge0$ which satisfies
the conditions of Theorem~1A we can write (by exploiting that
$\int f(x)\mu(\,dx)\le\frac\rho2\le\rho_k$) the estimate
$$
\mu_p(B_p(\Cal F))=N^{-p}_k |B_p(\Cal F)|\le N^{-p}_kV(\rho,p,N_0,k)
\le 2D\rho^{-p/4}
$$
by relation (2.8). Theorem~1A is proved.

\beginsection 3. The proof of Theorem 1.

First we prove the following Lemma 3.1 which is a special case of
Theorem 1.

\medskip\noindent
{\bf Lemma 3.1.} {\it Let us consider a finite set 
$X=\{x_1,\dots,x_{2^k}\}$ with $N=2^k$ elements together with
an $L_1$-dense class of function $\Cal F$ on $X$ with 
parameter $D\ge1$ and exponent $L\ge1$ that contains such 
functions $f\in\Cal F$ for which $0\le f(x)\le 1$ for all 
$x\in X$ and $\int f(x)\mu(\,dx)\le\rho$ with some 
$0<\rho<1$. Here $\mu$ denotes the uniform distribution 
on~$X$. Let us take the $n$-fold direct product $X^n$ of 
$X$ with some number $n\ge2$, and define the function 
$S_n(f)(x_{s_1},\dots,x_{s_n})=\summ_{j=1}^n f(x_{s_j})$
for all $(x_{s_1},\dots,x_{s_n})\in X^n$ and $f\in\Cal F$. 
Let us assume that $\rho\le n^{-200}$, and 
$N=2^k\ge\rho^{-3/2}$. Then the set $B_n(u)\subset X^n$ 
defined as
$$
B_n(u)=\left\{(x_{s_1},\dots,x_{s_n})\colon\; 
\sup_{f\in\Cal F}S_n(f)(x_{s_1},\dots,x_{s_n})>u\right\} \tag3.1
$$
satisfies the inequality
$$
\mu_n(B_n(u))\le 2D\rho^{u/25} \quad \text{for all }u\ge 40L, 
\tag3.2
$$
where $\mu_n$ denotes the uniform distribution on $X^n$.}

\medskip\noindent
{\it Proof of Lemma 3.1.}\/ Let us define for all functions 
$f\in\Cal F$  and integers $j$, $1\le j\le R$, where $R$ is 
defined by the relation $n<2^R\le 2n$, the functions 
$f_j(x)=\min(2^{-j},f(x))$ and $\bar f_j(x)=2^j f_j(x)$, 
$x\in X$. Put $\Cal F_j=\{f_j\colon\;f\in\Cal F\}$ and   
$\bar{\Cal F}_j=\{\bar f_j\colon\;f\in\Cal F\}$. One can simply 
check that $\Cal F_j$ is an $L_1$-dense class with parameter~$D$ 
and exponent~$L$, while $\bar{\Cal F}_j$ is an $L_1$-dense class 
with parameter~$D2^{jL}$ and exponent~$L$, if $\Cal F$ is 
an $L_1$-dense class with parameter~$D$ and exponent~$L$. We 
can also state that $\int f_j(x)\mu(\,dx)\le\rho$, and 
$\int \bar f_j(x)\mu(\,dx)\le 2^j\rho$ for all $f\in\Cal F$.

Let us define for all $f\in\Cal F$ and $1\le j\le R$ the 
following function $H_j(f)$ on $X^n$:
$$
H_j(f)(x_{s_1},\dots,x_{s_n})=\text{the number of such indices 
$l$ for which } \bar f_j(x_{s_l})=1. 
$$
We can write
$$
S_n(f)(x_{s_1},\dots,x_{s_n})\le\sum_{j=1}^R 
2^{1-j}H_j(f)(x_{s_1},\dots,x_{s_n})+1
$$
for all $f\in\Cal F$. This formula implies the inequality 
$$
\sup_{f\in\Cal F}S_n(f)(x_{s_1},\dots,x_{s_n})\le\sum_{j=1}^R 
2^{1-j}\sup_{f\in\Cal F}H_j(f)(x_{s_1},\dots,x_{s_n})+1,
$$
and the relation
$$
\align
&\left\{(x_{s_1},\dots,x_{s_n})\colon\; 
\sup_{f\in\Cal F} S_n(f)(x_{s_1},\dots,x_{s_n})>u\right\} \\
&\qquad \subset \bigcup_{j=1}^R\left\{(x_{s_1},\dots,x_{s_n})\colon\;
2^{1-j}\sup _{f\in\Cal F} H_j(f)(x_{s_1},\dots,x_{s_n})
>(\sqrt2-1)(u-1)2^{-j/2}\right\}.
\endalign
$$
Hence
$$
\mu_n(B_n(u))\le\sum_{j=1}^R \mu_n(D_n(u,j)) \tag3.3
$$
for the set $B_n(u)$ defined in~(3.1) by
$$
\align
D_n(u,j)=&\left\{(x_{s_1},\dots,x_{s_n})\colon\;\sup_{f\in\Cal F} 
H_j(f)(x_{s_1},\dots,x_{s_n})>\frac{\sqrt2-1}2(u-1)2^{j/2}\right\},\\
&\qquad\qquad\qquad\qquad\qquad,\qquad\qquad\qquad\qquad\qquad\qquad\qquad 
1\le j\le R.
\endalign 
$$
We can prove Lemma~3.1 with the help of relation~(3.3) if 
we give a good estimate on the measures $\mu_n(D_n(u))$. 
This can be done with the help of Theorem~1A.

Indeed, the set $D_n(u,j)$ consists of such sequences
$(x_{s_1},\dots,x_{s_n})\in X^n$ which have a subsequence
$(x_{s_{p_1}},\dots,x_{s_{p_t}})$ with 
$t=t(j)=[\frac{\sqrt2-1}2(u-1)2^{j/2}]+1$ elements, where 
$[\cdot]$ denotes integer part, with the property that 
there is a function $f\in\Cal F$ such that the function 
$\bar f_j(\cdot)$ defined with its help equals 1 in all 
coordinates of this subsequence. More explicitly,
$$
D_n(u,j)=\bigcup_{(\{l_1,\dots,l_t\}\subset\{1,\dots,n\}}
\left(\bigcup_{f\in\Cal F}
\{(x_1,\dots,x_n)\colon\;  \bar f_j(x_{s_{l_1}})=1,\dots,
\bar f_j(x_{s_{l_t}})=1\}\right) \tag3.4
$$
with $t=t(j)=[\frac{\sqrt2-1}2(u-1)2^{j/2}]+1$.

The outside union in (3.4) consists of 
$\binom n {t(j)}\le n^{t(j)}$ terms, and the cardinality 
of the sequences $(x_1,\dots,x_n)$ in the inner union 
can be bounded by means of Theorem~1A for each term if 
it is applied with $p=t(j)$, in the space $X$ consisting 
of $N=2^k=N_02^{\bar k}$ points, for the class of 
functions $\bar{\Cal F}_j$ which is an $L_1$-dense 
class of functions with parameter $D2^{jL}$ and 
exponent~$L$. Moreover, the functions 
$\bar f_j\in\bar{\Cal F}_j$ satisfy the inequality
$\int \bar f_j(x)\mu(\,dx)\le 2^j\rho$. This means that
under the conditions of Lemma~3.1 we can apply Theorem~1A
for the class of functions $\bar{\Cal F}_j$ with parameter
$\bar\rho=2^{j+1}\rho$ instead of~$\rho$. (We have to check
that all conditions of Theorem~1A hold. In particular, we 
can state that $\bar\rho=2^{j+1}\rho\le L^{-20}$, since 
$\rho\le n^{-200}$, $2^j\le 2n$, and since we estimate the 
probability in formula~3.2 only under the condition 
$u\ge 40L$, and this probability is zero if $u>n$, hence 
we may assume that $L\le\frac n{40}$. We chose the term 
$N_0$ in the application of Theorem~1A as $N_0=2^{k_0}$ 
with $k_0$ defined by the relation
$\frac1{16}\rho^{-3/2}< 2^{k_0}\le\frac18\rho^{-3/2}$, and 
$\bar k=k-k_0$.)

We will prove with the help of the above relations 
the inequality
$$
\aligned
\mu_n(D_n(u,j))&=\frac{|D_n(u,j)|}{N^n}
\le 2n^{t(j)} D2^{jL}(2^{j+1}\rho)^{t(j)/4}\\
&\le2D(8n^5\rho)^{t(j)/4}\le 2D\rho^{t(j)/5}\le D\rho^{ju/25}. 
\endaligned \tag3.5
$$
To get the first estimate in the second line of formula~(3.5) 
observe that under the condition of Lemma~3.1 
$\frac{\sqrt2-1}2(u-1)\ge 4L$, hence 
$2^{jL}\le 2^{j2^{-j/2}t(j)/4}\le2^{t(j)/4}$, and by the 
definition of the number $R$ we have  
$(2^{j+1})^{t(j)/4}\le (2^{R+1})^{t(j)/4}\le (4n)^{t(j)/4}$.
We imposed the condition $n\le \rho^{-1/200}$, and this 
implies the second inequality. Finally $t(j)\ge\frac{ju}5$. 
(In the last inequality a $j=1$ parameter is the worst case.) 
Relation (3.2) follows from~(3.3) and~(3.5). Lemma~3.1 is proved.

\medskip
Now we turn to the proof of the main result of this paper.

\medskip\noindent
{\it Proof of Theorem~1.}\/ We may assume that all functions 
$f\in\Cal F$ are non-negative, i.e. $0\le f(x)\le 1$ for all
 $f\in\Cal F$ and $x\in X$, because we can replace the function
$f$ by its absolute value $|f|$, and apply the result for
this new class of functions which also satisfies the conditions 
of Theorem~1. Next I show that we also may assume that the 
class of functions $\Cal F$ contains only finitely many
functions, satisfies the same conditions as the original
class of function $\Cal F$ with the only difference that
we assume that $\Cal F$ is an $L_1$-dense class with the same
exponent~$L$ but with parameter $D2^L$ instead of~$D$.

Indeed, if we have the same upper bound for the probability 
of $P\left(\supp_{f\in\Cal F'}S_n(f)>u\right)$ for all 
finite subsets $\Cal F'\subset \Cal F$, then this upper bound 
remains valid if we take the supremum for all $f\in\Cal F$. 
Besides, the conditions of Theorem~1 remain valid if 
$\Cal F$ is replaced by an arbitrary class of functions 
$\Cal F'\subset\Cal F$ with a small modification. Namely, 
we can state that $\Cal F'$ is an $L_1$-dense subclass 
with exponent $L$ but with a possibly different parameter 
$\bar D=D2^L$. (We had to change the parameter $D$ of an
$L_1$-dense class $\Cal F'\subset\Cal F$, because if a set 
of functions $f_1,\dots,f_m$ is an $\e$-dense class 
$\Cal F_{\e,\nu}$ appearing in the definition of $L_1$-dense 
property of the class of functions $\Cal F$, then these
functions $f_j$, $1\le j\le m$, may be not contained in 
$\Cal F'$. This problem can be overcome if we choose first 
an $\e/2$ dense subclass $\Cal F_{\e/2,\nu}$ in $\Cal F$ with
at most $D2^L\e^{-L}$ element, and then we replace the functions 
of this subclass with very close functions from $\Cal F'$ if
this is necessary.)

In the next step I show that we may restrict our 
attention to the case when the functions of the class of 
functions $\Cal F$ (consisting of finitely many functions)
take only finitely many values. For this goal first I split 
up the interval $[0,1]$ to $n$ subintervals of the following
form: $B_j=(\frac {j-1}n,\frac jn]$, $2\le j\le n$, and 
$B_1=[0,\frac1n]$. (We defined the function $B_1$
in a slightly different way in order to guarantee that 
the point zero is also contained in some set~$B_j$.) Then 
given a class of function $\Cal F$ on a set $X$ that contains 
finitely many functions $f_1,\dots,f_R$, we define the
following sets $A(s(1),\dots,s(R))\subset X$ (depending 
on~$\Cal F$):
$$
A(s(1),\dots,s(R))=\{x\colon\; f_j(x)\in B_{s(j)},\quad \text{for all }
1\le j\le R\},
$$
where $1\le s(j)\le n$ for all $1\le j\le R$. 

In such a way the sets $A(s(1),\dots,s(R))$ make up a partition 
of the set $X$. Actually, for the sake of a simpler argument 
we shall diminish a bit the set $X$, by defining it as the union 
of those sets $A(s(1),\dots,s(R))$ for which 
$\mu(A(s(1),\dots,s(R)))>0$ with the measure~$\mu$ appearing in
Theorem~1. This restriction will cause no problem in our 
later considerations. 

We shall define new functions $\tilde f_j(x)$, $1\le j\le R$,
by means of the partition of $X$ to the sets $A(s(1),\dots,s(R))$
by the formula
$$
\tilde f_j(x)=\frac {\int_{A(s(1),\dots,s(R))}f_j(x)\mu(\,dx)}
{\mu(A(s(1),\dots,s(R)))}, \quad 1\le j\le R, \quad \text{if }
x\in A(s(1),\dots,s(R)).
$$

We have $|f_j(x)-\tilde f_j(x)|\le\frac1n$ for all $1\le j\le n$
and $x\in X$. Hence 
$$
\left|\sup_{1\le j\le R} (S_n(f_j)
-S_n(\tilde f_j))\right|\le1, 
$$
for almost all sequences $\xi_1(\oo),\dots,\xi_n(\oo)$, and
as a consequence
$$
P\left(\supp_{1\le j\le R} S_n(f_j)>u+1\right)\le
P\left(\supp_{1\le j\le R} S_n(\tilde f_j)>u\right) \tag3.6
$$
Let us also observe that the class of functions 
$\tilde{\Cal F}=\{\tilde f_j,\;1\le j\le R\}$ also satisfies the 
conditions of Theorem 1, i.e. $\int \tilde f_j(x)\mu(\,dx)\le\rho$ 
for all $1\le j\le R$, and $\tilde{\Cal F}$ is an $L_1$-dense 
class with parameter~$\bar D=D2^L$ and exponent~$L$. (The conditions 
on the numbers $n$ and $\rho$ clearly remain valid.) 

The first relation follows from the identity 
$\int \tilde f_j(x)\mu(\,dx)=\int f_j(x)\mu(\,dx)$ which holds
because of the identities
$\int_{A(s(1),\dots,s(R))}\tilde f_j(x)\mu(\,dx)
=\int_{A(s(1),\dots,s(R))} f_j(x)\mu(\,dx)$ for all sets 
$A(s(1),\dots,s(R))$. 

To prove the $L_1$-dense property of $\tilde{\Cal F}$ let 
us introduce for all probability measures $\nu$ the probability 
measure $\tilde\nu=\tilde\nu(\nu)$ which is defined by the 
property that for all (measurable) sets $A(s(1),\dots,s(R))$ 
and $B\subset A(s(1),\dots,s(R))$ the identity
$\tilde\nu(B)=\mu(B)\frac{\nu(A(s(1),\dots,s(R))}
{\mu(A(s(1),\dots,s(R))}$ holds. Because of the special 
form of the  functions $\tilde f_j$ if a set of function 
$\tilde{\Cal F}_{\e,\tilde\nu}\subset\tilde{\Cal F}$ is an
$\varepsilon$-dense subset of $\tilde{\Cal F}$ in the space
$(X,\Cal X,\tilde\nu)$, then it is also $\varepsilon$-dense
in the space $(X,\Cal X,\nu)$. (In the proof of this 
statement we exploit that 
$$
\tilde\nu(A(s(1),\dots,s(R)))=\nu(A(s(1),\dots,s(R)))
$$
for all sets $A(s(1),\dots,s(R))$, and it depends only
on the value of a measure $\nu$ on the sets 
$A(s(1),\dots,s(R))$ whether a set of functions
$\{f_1,\dots,f_m\}\subset\tilde{\Cal F}$ is an $\e$-dense
subclass of $\tilde{\Cal F}$ with respect to the
measure $\nu$.)

Hence it is enough to prove the existence of an $L_1$-dense 
set $\tilde{\Cal F}_{\e,\nu'}$ with cardinality bounded by 
$\bar D\e^{-L}$ only with respect to such measures $\nu'$
which can be written in the form $\nu'=\tilde\nu(\nu)$ 
with some probability measure $\nu$. In this case the 
relation we want to check follows from the $L_1$-dense 
property of the original class of functions $\Cal F$ and 
the inequality 
$\int|\tilde f_j-\tilde f_{j'}|d\tilde\nu\le\int|f_j-f_{j'}|d\tilde\nu$
for all pairs $f_j,f_{j'}\in{\Cal F}_j$ 
and probability measure $\tilde\nu$. The last inequality holds, 
since
$$
\int_{A(s(1),\dots,s(R))}|\tilde f_j(x)-\tilde f_{j'}(x)|
\,d\tilde\nu(\nu(x)\le\int_{A(s(1),\dots,s(R))}|f_j(x)-f_{j'}(x)|
\,d\tilde\nu(x)
$$
for all sets $A(s(1),\dots,s(R))$.

Let us observe that for all $k\ge1$ we can define such a 
`discretized' probability measure $\bar\mu_k$ on the 
$\sigma$-algebra $\Cal X_k$ with atoms $A(s(1),\dots,s(R))$ 
in the space $X$ for which 
$$
|\bar\mu_k(A(s(1),\dots,s(R))-\mu(A(s(1),\dots,s(R))|\le 2^{-k},
$$
and
$$
\bar\mu_k(A(s(1),\dots,s(R))=\alpha(A(s(1),\dots,s(R)))2^{-k} \tag3.7
$$
with a non-negative integer $\alpha(A(s(1),\dots,s(R)))$
for all sets $A(s(1),\dots,s(R))$. (To find such a probability
measure $\bar\mu_k$ let us list the sets $A(s(1),\dots,s(R))$ 
as $B_1,\dots,B_Q$, and define the measure $\bar\mu_k$ by 
the relation $\summ_{l=1}^s \bar\mu_k(B_l)=\beta_s2^{-k}$ if 
$(\beta_s-1)2^k<\summ_{l=1}^s \mu_k(B_l)\le \beta_s2^{-k}$
with a positive integer~$\beta_s$. We assume this relation
for all $1\le s\le Q$.)

Clearly,
$$
P\left(\supp_{1\le j\le R} S_n(\tilde f_j)>u\right) 
=\lim_{k\to\infty}
P_{\bar\mu_k}\left(\supp_{1\le j\le R} S_n(\tilde f_j)>u\right)
\tag3.8
$$
for all $u>0$, where $P_{\bar\mu_k}$ means that we consider 
the probability of the same event as at the left-hand side of 
the identity, but this time we take iid. random variables 
$\xi_1,\dots,\xi_n$ with distribution $\bar\mu_k$ (on the 
$\sigma$-algebra generated by the atoms $A(s(1),\dots,s(R))$) 
in the definition of the random variables $S_n(\tilde f_j)$.

We shall bound the probabilities at the right-hand side in 
formula~(3.8) for all large indices~ $k$ by means of Lemma~3.1. 
This will be done with the help of the following construction. 
Take a space $\hat X=\hat X_k=\{x_1,x_2,\dots,x_{2^k}\}$ with 
$2^k$ elements and with the uniform distribution $\mu=\mu^{(k)}$ 
on its points. Let us fix a partition of $\hat X$ 
consisting of some sets $\hat A(s(1),\dots,s(R))$ with 
$\alpha(A(s(1),\dots,s(R)))$ elements, where the number 
$\alpha(\cdot)$ was introduced in~(3.7). Let us define the 
functions $\hat f_j(x)$, $1\le j\le R$, $x\in\hat X$, by 
the formula $\hat f_j(x)=\frac {s(j)}n$, $1\le j\le R$, 
if $x\in\hat A(s(1),\dots,s(R))$. Take the $n$-fold 
direct product $\hat X^n$ of $\hat X$ together with 
the uniform distribution $\mu_n=\mu_n^{(k)}$ on it and 
the functions $S_n(\hat f_j)(x_{t_1},\dots,x_{t_n})
=\summ_{l=1}^n \hat f_j(x_{t_l})$, $1\le j\le R$, 
if $(x_{t_1},\dots,x_{t_n})\in\hat X^n$ on the space
$\hat X^n$. I claim that
$$
\align
&P_{\bar\mu_k}\left(\supp_{1\le j\le R} S_n(\tilde f_j)>u\right)
\tag3.9 \\
&\qquad =\mu_n^{(k)}\left(\left\{(x_{t_1},\dots,x_{t_n})\colon\;
\supp_{1\le j\le R} S_n(\hat f_j)(x_{t_1},\dots,x_{t_n})>u\right\}\right)
\le 2\bar D\rho^{u/25} 
\endalign
$$
if $u>8L$.

The identity in formula (3.9) holds, since the joint distribution
of the random vectors $S_n(f_j)(\xi_1,\dots,\xi_n))$, $1\le j\le R$,
where $\xi_1,\dots,\xi_n$ are independent random variables  
with distribution $\bar\mu_k$ and of the random 
vectors
$S_n(\hat f_j)(x_{t_1},\dots,x_{t_n})$, $1\le j\le R$,
where the distribution of $(x_{t_1},\dots,x_{t_n})\in\hat X^n$ 
is $\mu_n^{(k)}$, agree. To prove the last inequality of~(3.9) 
it is enough to check that for all sufficiently large numbers~$k$ 
the class of functions 
$\hat{\Cal F}=\{\hat f_1,\dots,\hat f_R\}$ on the space 
$\hat X=\hat X_k$ satisfies the conditions of Lemma~3.1. Namely,
the $L_1$-dense property holds with parameter $\bar D=D2^L$ 
and exponent $L$, and $\int \hat f_j(x)\mu(\,dx)\le\rho$ 
with a number $\rho\le n^{-200}$ for all $\hat f_j\in\hat{\Cal F}$.

It is the $L_1$-dense property of the system $\hat X,\hat{\Cal F}$
that may demand some explanation. Let us observe that it is 
enough to check this property only for such probability 
measures $\hat\nu$ which have a constant density (with 
respect to the uniform distribution $\mu^{(k)}$) on all sets 
$\hat A(s(1),\dots,s(R))$. This reduction of the probability
measures can be justified similarly to the argument we applied
to prove the $L_1$-dense property of $\tilde\Cal F$ with the
help of the functions $\tilde\nu(\nu)$. Given a measure 
$\hat\nu$ on $\hat X$ with the above property let us 
correspond to it the measure $\tilde\nu$ on $X$ defined by 
$\tilde\nu(A(s(1),\dots,s(R))=\hat\nu(\hat A(s(1),\dots,s(R))$
for all sets $A(s(1),\dots,s(R))$. Then we get that if a 
class of functions $\tilde{\Cal F}_{\e,\tilde\nu}
=\{\tilde f_{l_1},\dots,\tilde f_{l_s}\}$ is an
is an $\varepsilon$-dense class of $\tilde{\Cal F}$
with respect to the measure $\tilde\nu$, then the class
of function $\hat{\Cal F}_{\e,\hat\nu}
=\{\hat f_{l_1},\dots,\hat f_{l_s}\}$ is an is an 
$\varepsilon$-dense class with respect to the 
measure~$\hat\nu$. The $L_1$-density property of 
$\hat\Cal F$ follows from this fact.

Then we get the inequality part of formula~(3.9) from 
Lemma~3.1. Relation~(1.1) follows from~(3.9), (3.8) and~(3.6). 
We still have to understand that in our estimation the 
coefficient $2\bar D=2D2^L$ in~(3.9) can be replaced by~$D$ if 
we estimate the probability~(1.1) only for $u\ge \frac14(L+1)$, 
and the term $\rho^{u/25}$ in~(3.9) is replaced by 
$\rho^{u/50}$ when turning from~(3.9) to formula~(1.1). 
To see this observe that 
$\rho^{u/25}\le\rho^{\frac14(L+1)/50}\cdot\rho^{u/50}\le 
n^{-(L+1)}\rho^{u/50}\le \frac12 2^{-L}\rho^{u/50}$ if $u\ge\frac14L$, 
$\rho\le n^{-200}$, and $n\ge2$. Theorem~1 is proved.

\vfill\eject

\noindent
{\bf References.}

\medskip
\item{[1]} J. Koml\'os, P. Major, G. Tusn\'ady, An approximation
of partial sums of independent rv.'s and the sample DF. II
Z. Wahrscheinlichkeitstheorie verw. Gebiete 34, 33--58 (1976)
\item{[2]} P. Major  On the tail behaviour of the distribution
function of the supremum of a class of partial sums of i.i.d. 
random variables. submitted to Electron. J. of Probab.
\item{[3]} D. Pollard, {\it Convergence of Stochastic Processes}
(Springer, New York, 1984)

\bye